\documentclass[notitlepage,11pt]{article}
\usepackage{amssymb,amsfonts,amsmath,geometry,esint,comment}% upgreek,cases,graphicx,ulem}
\usepackage{scalerel,stackengine,mathtools,stmaryrd,color}

\usepackage{color}

\catcode`\@=11
\@addtoreset{equation}{section}

\catcode`\@=12

\def\eps{\varepsilon}

\def\R{\mathbb{R}}

\def\f{\varphi}
\def\e{\varepsilon}

\def\D{\mathcal D}

\def\S{\mathbb S}
\def\Hper{H^1_{\rm per}}
\def\Hdual{H^{-1}_{\rm per}}

\def\p{\overline u}
\def\irn{\fint\limits_{\S^1}}
\def\irndue{\int\limits_{\R^2}}

\def\D{\mathbb D}
\def\cmp{{\rm c}_{mp}}

\DeclareMathOperator*{\esssup}{ess\,sup}

\def\div{{\rm div}}

\def\proof{\noindent{\textbf{Proof. }}}
\def\QED{\hfill {$\square$}\goodbreak \medskip}

\newtheorem{Theorem}{Theorem}[section]
\newtheorem{Lemma}[Theorem]{Lemma}
\newtheorem{Proposition}[Theorem]{Proposition}
\newtheorem{Corollary}[Theorem]{Corollary}

\newtheorem{Definition}[Theorem]{Definition}

\newtheorem{Claim}{Claim}

\def\ddiv{{\rm d\!l}}

\def\rmH{{\rm H}}

\begin{document}

\title{Planar loops with prescribed curvature via Hardy's inequality}

\author{{Gabriele Cora\footnote{Dipartimento di Matematica "G. Peano", Universit\`a di Torino, Italy. Email: {gabriele.cora@unito.it}, orcid.org/0000-0002-0090-5470.}}, 
Roberta Musina\footnote{Dipartimento di Scienze Matematiche, Informatiche e Fisiche, Universit\`a di Udine, Italy.
Email: {roberta.musina@uniud.it}, orcid.org/0000-0003-4835-8004.}
}

\date{}

\maketitle

\noindent
{\small {\bf Abstract:} We investigate the existence of closed planar loops with prescribed curvature. Our approach is variational, and relies on a Hardy type inequality and its associated functional space.}



\section{Introduction}

Let $\rmH:\R^2\to \R$ be a given  function. We study the existence of {\em $\rmH$-loops}, which are 
solutions to
\begin{equation}
\label{eq:problem}
\begin{cases}
u''=|u'| \rmH(u)iu' \\
u\in C^2(\S^1,\R^2)~,~~\textit{$u$ non constant.}
\end{cases}
\end{equation}
If $u$ solves \eqref{eq:problem}, then $u'$ is orthogonal
to $u''$. Hence $|u'|$ is constant and $u$ is a regular curve, having curvature $\rmH(u)$ at each point.

The problem of the existence of $\rmH$-loops has been raised in \cite[Question $(Q_0)$]{BCG}. Our interest in (\ref{eq:problem}) is also motivated by its relation with Arnold's problem on magnetic geodesics 
\cite[Problems 1988/30, 1994/14 and 1996/18]{A2}.

It is easy to see that (\ref{eq:problem}) has no solution if $\rmH=0$; 
if $\rmH\neq 0$ is constant then $u$  solves (\ref{eq:problem})
if and only if $u$ parametrizes a circle of radius $1/|\rmH|$ with constant scalar speed $|u'|$.
%Despite 
In contrast and despite its simple formulation, the variable curvature case is more involved. Indeed, nonexistence phenomena may occur \cite{KL} and only a few  results are available in the literature, see \cite{CalCo, CG, GR, MuLoop, MuZ} and references therein.

The main novelty in our approach consists in making the most of Hardy's inequality, via the introduction of the quantity
\begin{equation}
\label{eq:NH}
N_\rmH:=\frac{1}{\sqrt{4\pi}}\Big(\irndue|\nabla \rmH(z)\cdot z|^2~\!dz\Big)^\frac12.
\end{equation}

If $\rmH$ is differentiable, then $\nabla \rmH(z)\cdot z$ only depends on $|z|$ and on the radial derivative of $\rmH$ at $z$;
otherwise, the expression $\nabla \rmH(z)\cdot z$ in (\ref{eq:NH}) has to be interpreted in a distributional (or weak) sense.
We refer to Section \ref{S:ddiv}, and in particular to Definition \ref{D:def_ddiv}, for details and related observations.

We use variational methods to prove the next existence result.

\begin{Theorem}
\label{T:main}
Assume that $\rmH\in C^0(\R^2)$ satisfies

\smallskip

$(H_1)$~~$N_\rmH<1$,

\smallskip

$(H_2)$~~$\rmH(z)-1=o(|z|^{-1})$ as $|z|\to\infty$,

\smallskip

$(H_3)$~~there exist $\tilde{p} \in \R^2$ such that 
$\rmH(z)\ge  1$ if $|z-\tilde{p}|< 2(1 + N_\rmH)$.

\smallskip

\noindent
Then there exists at least one $\rmH$-loop.
\end{Theorem}

Let us comment our hypotheses in relation to the existing literature. 

The variational approach to problem (\ref{eq:problem}) is well understood, see for instance \cite[Section 1.3]{BCG}. The  energy functional takes the form
$$
E_\rmH(u)=L(u)+A_\rmH(u)~,\qquad u\in H^1(\S^1,\R^2)~\!.
$$
Here 
$$
L(u):=\Big(\irn|u'|^2~\!d\theta\Big)^\frac12\ge \frac{1}{2\pi}\cdot(\textit{lenght of $u$})~\!,
$$
and $A_\rmH(u)$ is proportional to the  {\em $\rmH$-area functional}. It measures the algebraic area enclosed by $u$ with respect to the weight $\rmH$.
Details on $A_\rmH$ can be found in Section \ref{S:area}. 

The functional $E_\rmH$ is Fr\'echet differentiable on 
$H^1(\S^1,\R^2)\setminus\R^2$, and
any critical point $u$ of $E_\rmH$ is a weak solution to (\ref{eq:problem}). 
It is then easy to check that $u$ is in fact 
an $\rm H$-loop.

Reasonable hypotheses on $\rmH$ allow to construct candidate critical levels for $E_\rmH$.
However, severe lack of compactness phenomena could occur. The worst ones are produced by the group of dilations in the target space $\R^2$:
in Claim \ref{E:Ex2} of the Appendix, we show that there might exist Palais-Smale sequences $u_n$ having unbounded seminorms $L(u_n)$. Even more impressive phenomena 
have been observed in the related and more challenging {\em $\rmH$-bubble problem} (see the collection of examples in \cite{CM_Adv}), due to 
the interaction between the groups of dilations in the target space and of M\"obius transforms in the domain.

In \cite{MuLoop}, as well as in the papers \cite{Paolo_JFA, CM_CCM, CM_Arch, CM_NATMA} on the $\rmH$-bubble problem, this loss of compactness was addressed by imposing $\rmH\in C^1(\R^2)$ and $M_\rmH<1$, where
\begin{equation}
\label{eq:MH}
M_\rmH:=\sup_{z\in\R^2}|(\nabla \rmH(z)\cdot z)z|~\!.
\end{equation}
The assumption  $M_\rmH<1$ also affects the topology of the energy sublevels and  the properties of the Nehari manifold
$\Sigma=\{dE_\rmH(u)u=0\}$. In particular, $\Sigma$ turns out to be a smooth and natural constraint for $E_\rmH$,
so that any minimizer for $E_\rmH$ over $\Sigma$ gives rise to an $\rm H$-loop of minimal energy. In \cite[Theorem 2.5]{MuLoop}, the existence of
a minimal $\rmH$-loop is obtained by asking, in addition,  that  
$(H_2)$ in Theorem \ref{T:main} is satisfied, together with 
\begin{itemize}
\item[$(*)$] there exist $R>0$  such that 
$\rmH(z)\ge  1$ for $|z|>R$.
\end{itemize}
Assumption $(*)$, which is evidently stronger than $(H_3)$, is asked in \cite{MuLoop} 
to prevent the lack of compactness produced by the group of traslations in $\R^2$.
It is important to stress that the smoothness assumption $\rmH\in C^1(\R^2)$ cannot be easily removed via an approximation argument.

\medskip

Our starting goal was to find 
an alternative to the "$L^\infty$-type" hypothesis $M_\rmH<1$ which, among others, would allow to include non differentiable curvatures. 
We succeeded in this purpose by introducing the constant $N_\rmH$ and the assumption $(H_1)$.
It is important to notice that the hypotheses $M_\rmH<1$ and $N_\rmH<1$ are not comparable even in the case $\rmH\in C^1(\R^2)$, see
Claim \ref{E:Ex1} in the Appendix.

We sketch here the main steps in the proof of Theorem \ref{T:main}.

First, we show that any continuous non constant curvature satisfying $N_\rmH<\infty$ and $(H_2)$ obeys the Hardy type inequality
\begin{equation}
\label{eq:Hardy0}
\irndue| \rmH(z)-1|^2~\!dz\lneq \irndue|\nabla \rmH(z)\cdot z|^2~\!dz~\!.
\end{equation}
This is proved in Section \ref{S:ddiv} via a density result (see Lemma \ref{L:dense}), of independent interest.
We then point out,  in Section \ref{S:area}, some noteworthy properties of the area functional $A_\rmH$.

Next, in Lemma \ref{L:HPS} we show that if in addition $N_\rmH <1$, then the Palais-Smale condition fails
only at  energy levels $\ell/2$, $\ell\in\mathbb N$. This is the main step in the proof of Theorem \ref{T:main},  to which
Section \ref{S:energy} is dedicated.

Since $\rmH$ is not required to be differentiable, then the Nehari manifold is not smooth (and is not a natural constraint),
so that the approach used in \cite{MuLoop} fails. However, it is possible to construct 
in an almost standard way a positive mountain-pass energy level $\cmp$, compare with (\ref{eq:mp}).

The assumption $(H_3)$ is needed only in the last step of the proof, to show that either  
there exists a circle of radius $1$ which can be parametrized by an $\rm H$-loop, or $\cmp<\frac12$. In the latter
case, we have enough compactness to infer the existence of an $\rm H$-loop, which is a mountain-pass critical point for the energy functional  $E_\rmH$.

\medskip

Finally, we notice that the assumption $(H_2)$ can be replaced by $\rmH(z)-\lambda
=o(|z|^{-1})$ for some constant $\lambda \neq 0$, after suitably modifying $(H_3)$; in contrast, it turns
out that problem (\ref{eq:problem}) has no solutions if $\lambda=0$, see Corollary
\ref{C:Hinfty} and Theorem \ref{T:ne}, respectively.

{\small
\paragraph{Notation.}
The Euclidean space $\R^2$ is endowed with the scalar product
$z\cdot \xi$ and norm $|\cdot|$. The standard $L^p$-norm, $p\in[1,\infty]$, is denoted by $\|\cdot\|_p$.

We denote by $\D_r(z)$ the disc in $\R^2$ of radius $r >0$  about $z \in \R^2$. We simply write $\D_r$
instead of $\D_r(0)$.

We will often use complex notation for points in $\R^2$. For instance, we write $iz=(-y,x)$ for $z=(x,y)\in\R^2$ and 
put
$\S^1:=\partial\D_1\equiv \{e^{i\theta}~|~\theta\in\R\}$.

If $f$ is a differentiable function on $\S^1$, we put $f'(\sigma)=f'(\sigma)(i\sigma)$, so that $f'$ is a function on $\S^1$ as well.
In fact, any function $f$ on $\S^1$ can be identify with a $2\pi$-periodic function via the identity
$f(\theta)\equiv f({e^{i\theta}})$. For this reason, we put
$$
{\Hper}:=H^1(\S^1,\R^2).
$$
We identify constant functions $\S^1\to \R^2$ and  points in $\R^2$. Thus $\Hper\setminus\R^2=\{L(u)\neq 0\}$ contains only non constant functions.
We endow $\Hper$ with the equivalent norm
\begin{equation}
\label{eq:mean}
\|u\|_{H^1_{\rm per}}^2 = L(u)^2 + |\overline u |^2  \,,\quad \text{where}\quad 
L(u) = \Big(\irn|u'|^2\, d\theta\Big)^2\,, \quad \overline u = \irn u\, d \theta\,,
\end{equation}
and denote by $\Hdual$ its dual space. 
Recall that $H^1_{per} \hookrightarrow C^0(\S^1,\R^2)$ with compact embedding.
For future convenience we point out the elementary (and rough) inclusion 
\begin{equation}
\label{eq:Linfty}
u(\S^1)\subset \D_{\rho_{\!u}}(\p)~,\quad \text{where $\rho_u:=2\pi L(u)$.}
\end{equation}}

\section{The Hardy inequality and the quantity $N_\rmH$} 
\label{S:ddiv}

Our approach is crucially based on the classical {sharp} Hardy inequality
\begin{equation}
\label{eq:Hardy}
\irndue|K(z)|^2~\!dz \lneq \irndue |\nabla K(z)\cdot z|^2~\!dz,
\end{equation}
which holds for any nontrivial function $K\in C^\infty_c(\R^2)$.  
The first order differential operator in (\ref{eq:Hardy}) and its formal 
adjoint, 
$$
\ddiv K(z):=\nabla K(z)\cdot z~,\qquad \ddiv^* K(z):=-\div(K(z)z)~\!,
$$
can be extended in a standard way
to ${K}\in L^1_{\rm loc}(\R^2)$ in the sense of distributions. More precisely, for 
 $\f\in C^\infty_c(\R^2)$ we define
$$
\begin{gathered}
\langle{\ddiv K},\f\rangle:=-\irndue {K}(z)\div(\f(z)z)~\!dz
=\langle{K},\ddiv^*\f\rangle~,\quad
\langle\ddiv^* {K},\f\rangle:=\irndue {K}(z)(\nabla \f(z)\cdot z)~\!dz = \langle{K},\ddiv\f\rangle.
\end{gathered}
$$

We can now give a precise explanation of the hypothesis $(H_1)$ in the introduction.
\begin{Definition}
\label{D:def_ddiv}
Let $\rmH\in L^1_{\rm loc}(\R^2)$. We say that $\ddiv\rmH\in L^2(\R^2)$ if the distribution $\ddiv\rmH$
can be continuously extended to $L^2(\R^2)$. In this case, we put
$$
\irndue|\nabla \rmH(z)\cdot z|^2~\!dz:=\irndue|\ddiv\rmH(z)|^2~\!dz.
$$
\end{Definition}

In order to prove Theorem \ref{T:main} we need to show that any nonconstant curvature $\rmH$ satisfying the assumptions 
therein obeys Hardy's inequality (\ref{eq:Hardy0}). More generally,
the next lemma holds.

\begin{Lemma}\label{L:HinD_new}
Let $\rmH \in C^0(\R^2)$ be non constant. Assume that $N_\rmH<\infty$ and that $(H_2)$ is satisfied. Then  $\rmH-1\in L^2(\R^2)$, $\ddiv\rmH \in L^2(\R^2)$ and
the Hardy type inequality (\ref{eq:Hardy0}) holds.
\end{Lemma}

Lemma \ref{L:HinD_new} is in fact an immediate corollary of Lemmata \ref{L:HinD} and \ref{L:dense} below. Their proofs require a preliminary result, which could be of independent interest.

\begin{Lemma}
\label{L:dense2}
Let $K\in L^2_{\rm loc}(\R^2)$ be such that $\ddiv K\in L^2_{\rm loc}(\R^2)$.
\begin{itemize}
\item[$i)$] There exists a sequence $K_\eps\in C^\infty(\R^2)$ 
such that $K_\eps\to K$ and $\ddiv K_\eps\to \ddiv K$ in $L^2_{\rm loc}(\R^2)$;
\item[$ii)$] $\ddiv(\psi K)=\psi\ddiv K+K\ddiv\psi$ for any $\psi\in C^\infty_c(\R^2)$.
\end{itemize}
\end{Lemma}

\proof
Let $(\rho_\eps)_\eps$ be a sequence of radially decreasing 
mollifiers.  Then 
$K_\eps:=K*\rho_\eps\in C^\infty(\R^2)$ and $K_\eps\to K$ in $L^2_{\rm loc}(\R^2)$. 
To prove $i)$ it remains to show that $\ddiv K_\eps\to \ddiv K$ in $L^2_{\rm loc}(\R^2)$.

For any $z\in\R^2$ we compute 
$$
\begin{aligned}
\big({K}*(\ddiv\rho_\eps)\big)(z)
&=\irndue {K}(\xi)\nabla\rho_\eps(z-\xi)\cdot(z-\xi)~\!d\xi\\
&=\irndue {K}(\xi)\nabla\rho_\eps(z-\xi)\cdot z~\!d\xi
-\irndue {K}(\xi)\nabla\rho_\eps(z-\xi)\cdot\xi~\!d\xi
\\
&=\nabla({K}*\rho_\eps)(z)\cdot z
+((\ddiv^*{K})*\rho_\eps)(z)
=\ddiv({K}*\rho_\eps)(z)
+((\ddiv^*{K})*\rho_\eps)(z)
\end{aligned}
$$
to infer the identity
\begin{equation}
\label{eq:identity}
\ddiv K_\eps=-(\ddiv^*{K})*\rho_\eps+{K}*(\ddiv\rho_\eps).
\end{equation}
Firstly, we get that $-(\ddiv^*{K})*\rho_\eps=-\ddiv^*{K}+o(1)$ in $L^2_{\rm loc}(\R^2)$, as $-\ddiv^*{K}=\ddiv{K}+2{K}\in L^2_{\rm loc}(\R^2)$. 
Then we use integration by parts to get 
$$
-\frac12\irndue \ddiv\rho_\eps(z)~\!dz=-\frac12\irndue \nabla\rho_\eps(z)\cdot z~\!dz =1.
$$
We see that $\big(-\frac12\ddiv\rho_\eps)_\eps$
is a sequence of (nonnegative) mollifiers as well (recall that $\rho_\eps$ is radially decreasing), hence ${K}*(\ddiv\rho_\eps)=-2{K}+o(1)$ in $L^2_{\rm loc}(\R^2)$.
In conclusion,  from (\ref{eq:identity}) we obtain that 
$\ddiv K_\eps=-\ddiv^*{K}-2{K}+o(1)=\ddiv{K}+o(1)$ in $L^2_{\rm loc}(\R^2)$, which ends the proof of $i)$.

\medskip
Finally, let $K_\e\in C^\infty(\R^2)$ be the sequence in $i)$. Since trivially
$$
\ddiv(\psi K_\e)=\psi\ddiv K_\e+K_\e\ddiv\psi~,\quad \text{for any $\psi\in C^\infty_c(\R^2)$,}
$$
then $ii)$ follows by taking the limit as $\eps\to 0$.
\QED

For completeness we include below, in addition to the statements which are needed to prove Theorem \ref{T:main}, some side results and observations.

We introduce the domain of the unbounded and densely defined operator $\ddiv$ on $L^2(\R^2)$, namely
$$
\widehat{\mathcal D}^1(\R^2)=\big\{{K}\in L^2(\R^2)~|~\ddiv {K}\in L^2(\R^2)\big\}\,.
$$
Clearly, $\widehat{\mathcal D}^1(\R^2)$ is a Hilbert space with respect to the norm
$\|K\|^2=\|\ddiv {K}\|_2^2+\|K\|_2^2$. Notice that $\widehat{\mathcal D}^1(\R^2)$ is larger than the 
nowadays standard weighted homogeneous space
$$
\mathcal D^1(\R^2;|z|^2 dz)=\big\{ {K}\in L^2(\R^2)~|~|\nabla {K}|\in L^2(\R^2;|z|^2dz)\big\}. 
$$
Evidently, a radial function belongs to $\widehat{\mathcal D}^1(\R^2)$ if and only if it belongs to $\mathcal D^1(\R^2;|z|^2 dz)$. It is well known that (\ref{eq:Hardy}) holds with a sharp constant for any  $K\in \mathcal D^1(\R^2;|z|^2 dz)\setminus\{0\}$.

\begin{Lemma}
\label{L:dense}
The following facts hold.
\begin{itemize}
\item[$i)$] $C^\infty_c(\R^2)$ is dense in $\widehat{\mathcal D}^{1}(\R^2)$;
\item[$ii)$] $\|\ddiv K\|_2=\|\ddiv^*K\|_2$ for any $K\in \widehat{\mathcal D}^1(\R^2)$;
\item[$iii)$] The Hardy type inequality (\ref{eq:Hardy}) holds with a sharp constant. 
Therefore, $\|\ddiv \cdot \|_2$ is an equivalent Hilbertian norm on $\widehat {\mathcal D}^1(\R^2)$.
\end{itemize}
\end{Lemma}

\proof
By exploiting the proof of $i)$ in Lemma \ref{L:dense2} one can show that $C^\infty(\R^2)\cap \widehat{\mathcal D}^{1}(\R^2)$
is dense in $\widehat{\mathcal D}^{1}(\R^2)$. Then $i)$ follows in a standard way. Details
are omitted.

\medskip

Next, let $K\in C^\infty_c(\R^2)$. We use integration by parts to compute 
$$
\irndue |\nabla K(z)\cdot z|^2~\!dz= \irndue |\div( K(z) z)-2K(z)|^2~\!dz= \irndue |\div( K(z) z)|^2~\!dz.
$$
To conclude the proof of $ii)$ use the density result in $i)$.

Next, a standard way to prove the classical Hardy inequality goes as follows. 
Use integration by parts to get 
$$
\irndue K\ddiv K~\!dz=
-\irndue |K|^2~\!dz~
$$
for any $K\in C^\infty_c(\R^2)$. Thus $\| K\|^2_2\le \|\ddiv K\|^2_2$ for any $K\in \widehat{\mathcal D}^{1}(\R^2)$ by the density result in $i)$. Therefore, by known facts we have that (\ref{eq:Hardy}) holds with a sharp constant and with a large inequality instead of the strict one. 

Assume that $K_0\in \widehat{\mathcal D}^{1}(\R^2)$ satisfies $\|\ddiv K_0\|_2=\|K_0\|_2$.  Since it achieves the minimum of the map $K\in \widehat{\mathcal D}^{1}(\R^2)\mapsto \|\ddiv K\|^2_2-\|K\|^2_2$, then it solves
\begin{equation}
\label{eq:EL}
\ddiv^*\ddiv K=  K
\end{equation}
in a weak sense. Let $\overline K$ be the $L^2$-orthogonal projection of $K$ on the space of radial functions.  
Then $ \overline K\in \mathcal D^1(\R^2;|z|^2~\!dz)$.
Testing (\ref{eq:EL}) with $\overline K$, we see that $\|\ddiv^*\overline K\|_2=\|\overline K\|_2$. This
implies that $K_0\equiv 0$,  as the Hardy constant is not achieved on $\mathcal D^1(\R^2;|z|^2~\!dz)$, and
completes the proof.
\QED

The last result in this section, together with Lemma \ref{L:dense}, readily implies Lemma \ref{L:HinD_new}, by choosing $K=\rmH-1$.

\begin{Lemma}\label{L:HinD}
Let $K\in L^2_{\rm loc}(\R^2)$. Assume that $\ddiv K\in L^2(\R^2)$, and that the function $z \to |z|K(z)$ is in $L^\infty(\R^2)$. Then $K\in \widehat{\mathcal D}^1(\R^2)$.
\end{Lemma}

\proof
Take a cut-off function $\psi\in C^\infty_c(\R^2)$ such that $0\le \psi\le 1, \psi\equiv 1$ on $\{|z|< 1\}$ and $\psi\equiv 0$ on
$\{|z|>2\}$. For any integer $h\ge 1$ put $\psi_h(z)=\psi\big(\frac{z}{h}\big)$. Evidently 
$\psi_h K\to K$ in $L^2_{\rm loc}(\R^2)$.

By $ii)$ in Lemma \ref{L:dense2} we have that 
 $\psi_h K\in \widehat{\mathcal D}^1(\R^2)$ and
 $$
\|\ddiv(\psi_hK)\|_2=\|\psi_h\ddiv K+K\ddiv\psi_h\|_2\le \|\ddiv K\|_2+ \|K\ddiv\psi_h\|_2.
$$
Note that $|\ddiv\psi_h(z)|=|\nabla \psi_h(z)\cdot z|\le 2\|\nabla\psi\|_\infty$ on $\text{supp}(\psi_h)\subseteq\{h\le|z|\le 2h\}$. Thus
$$
\|K\ddiv\psi_h\|_2^2\le c_\psi \int\limits_{\{h<|z|<2h\}}|K|^2~\!dz\le
c_\psi \int\limits_{\{h<|z|<2h\}}|zK(z)|^2|z|^{-2}~\!dz\le c_\psi \esssup_{z \in \R^2}|zK(z)|^2\,,
$$
where the constants $c_\psi$ depend only on $\psi$.
We showed that the sequence $(\psi_h K)_h$ is bounded in $\widehat{\mathcal D}^1(\R^2)$ endowed with the norm $\|\ddiv \cdot \|_2$ (compare with $iii)$ in Lemma \ref{L:dense}), which is sufficient to conclude that 
$K\in \widehat{\mathcal D}^1(\R^2)$. 
\QED

\section{Area functionals}
\label{S:area}

In this section we collect some partially known results about  {\em $K$-weighted area functionals}.

\medskip
Assume firstly that 
$K$ is a given function in $C^\infty(\R^2)$  and take a vectorfield $Q\in C^1(\R^2,\R^2)$ such that $\div Q=K$. For instance, choose
$Q(z)=z\displaystyle{\int_0^1 K(sz)s~\!ds}$.

\medskip
The functional
\begin{equation}
\label{eq:AKu}
A_K(u)=\irn Q(u)\cdot iu'~\!d\theta~,\quad u\in\Hper
\end{equation}
is well defined and weakly continuous on $\Hper$, use the compactness of the embedding
of $\Hper$ into $C^0(\S^1,\R^2)$. Moreover, it is of class $C^1$ on $\Hper$, with differential given by
\begin{equation}
\label{eq:dA}
dA_K(u)=\frac{1}{2\pi}K(u)iu'.
\end{equation}
For the proof, use integration by parts and the  identity
$
(dQ(u)\f)\cdot iu' -(dQ(u)u')\cdot i\f=\div Q(u)~\!\f\cdot iu'$.

Since $A_K$ vanishes on constant functions, we see
that $A_K(u)$ does not depend on the choice of $Q$.

\medskip

If $K\equiv 1$ then Fourier series can be used to prove that the $1$-area functional 
$$
{A_1}(u)=\frac12\irn u\cdot iu'~\!d\theta
$$
is analytic on $\Hper$, and satisfies
\begin{equation}
\label{eq:iso}
2|{A_1}(u)|\le L(u)^2~,\quad 
dA_1(u)= \frac{1}{2\pi}iu'~,\quad dA_1(u)u=2A_1(u)\qquad 
\text{for any $u\in {\Hper}$.}
\end{equation}

If $K\in C^\infty_c(\R^2)$ we can take $Q=\nabla V_K$, where $V_K=\frac{1}{2\pi}(K*\log|\cdot|)$ is the solution to the Poisson equation 
$\Delta V_K=K$. It readily follows that
\begin{equation}
\label{eq:div}
{A}_K(u)=\frac{1}{2\pi}\int\limits_{\R^2}K(z)j_u(z)~\!dz\quad \text{for any $u\in \Hper$,}
\end{equation}
where 
$$
j_u(z)=\irn\frac{u-z}{|u-z|^2}\cdot iu'~\!d\theta\qquad \text{for $z\notin u(\S^1)$}
$$
is the winding number of the loop $u-z$. Trivially, $j_u$ takes only integer values and vanishes outside 
any disk  containing $u(\S^1)$. In particular, by (\ref{eq:mean}) and (\ref{eq:Linfty}) we have
\begin{equation}
\label{eq:j_supp}
\text{supp}(j_u)\subset  \D_{\rho_u}(\overline u)~.
\end{equation}

The estimate on $j_u$ in the next proposition is crucially used in our approach.

\begin{Proposition}
\label{P:u_fixed}
Let $u\in {\Hper}$. Then $j_u\in L^2(\R^2,\mathbb Z)$ and
$$
\|j_u\|_2\le \sqrt{\pi}\irn|u'|~\!d\theta\le \sqrt{\pi} L(u).
$$
\end{Proposition}

\proof
We provide an alternative to the proof of Theorem 3 in \cite{CC}. Our argument 
has been inspired by 
\cite[Section 2]{Ste} and is based on 
Federer's theory of integral currents \cite{Fe}.

We identify functions ${K}\in C^\infty_c(\R^2)$ with $2$-forms ${K}(z)dx\wedge dy$ on $\R^2$, 
and introduce the $2$-dimensional current $J_u$ given by 
$$
J_u(K):=\frac{1}{2\pi}\irndue j_u(z) K(z)dx\wedge dy~\!.
$$ 
Notice that (\ref{eq:div}) becomes
$$
J_u(K)=A_K(u).
$$
The boundary of $J_u$ is 
the $1$-current defined via
$\partial J_u(\alpha)=J_u(d\alpha)$. 
Given a $1$-form $\alpha$, we take the vectorfield $Q^\alpha$ such that $\alpha=-Q^\alpha_2(z)dx+Q^\alpha_1(z)dy$. Thus $d\alpha=(\div Q^\alpha)dx\wedge dy$ and
$$
\partial J_u(\alpha)=J_u(\div Q^\alpha)=A_{\div Q^\alpha}(u)=
\irn Q^\alpha(u)\cdot iu'~\!d\theta. 
$$
This allows us to estimate the mass of $\partial J_u$ by 
$$
M(\partial J_u(\alpha))=
\sup_{\|\alpha\|_\infty=1}|\partial J_u(\alpha)|\le \sup_{\|Q\|_\infty=1}~\!\big|\irn Q(u)\cdot iu'~\!d\theta\big|
\le  \irn|u'|~\!d\theta.$$ 
Since $j_u$ has compact support, then the
conclusion follows by  \cite[Theorem 4.5.9, statement (31)]{Fe}.
\QED

We summarize in the next lemmata few consequences of the previous observations. The first one
readily follows from  formulae (\ref{eq:div}), (\ref{eq:j_supp}) and  Proposition \ref{P:u_fixed}, thanks to the density 
of $C^\infty_c(\R^2)$ in $L^2(\R^2)$.

\begin{Lemma} 
\label{L:2.5}
Let $u\in {\Hper}$. Then the following facts hold. 
\begin{itemize} 
\item[$i)$] The area functional $A_K(u)$ in (\ref{eq:div}) is well defined for any $K\in L^2_{\textrm{loc}}(\R^2)$ and
the linear map $K\to A_K(u)$ can be continuously extended to $L^2(\R^2)$;
\item[$ii)$]  If $K\in L^2(\R^2)$, then 
$$\sqrt{4\pi}|A_K(u)|\le \big(\hskip-0.2cm\int\limits_{\D_{\rho_u}\!(\p)}\hskip-0.2cm|K(z)|^2~\!dz\big)^\frac12~\!L(u),$$
where $\rho_u=2\pi L(u)$.
In particular, the following weighted isoperimetric inequality holds,
\begin{equation}
\label{eq:isop}
\sqrt{4\pi}|A_{K}(u)|\le   \|{K}\|_{2} L(u);
\end{equation}
\end{itemize}
\end{Lemma}

\begin{Lemma} 
\label{L:2.6}
Let $K\in C^0(\R^2)$. Then the functional $u\mapsto A_K(u)$ is continuously differentiable on $\Hper$ and formula (\ref{eq:dA}) holds.

If in addition $K\in \widehat{\mathcal D}^1(\R^2)$, then 
$dA_K(u)u = A_{\ddiv K}(u) + 2 A_K(u)=-A_{\ddiv^*\!K}(u)$ for any $u\in \Hper$.
\end{Lemma}

\proof 
Let  $(\rho_\eps)_\eps$ be a sequence of mollifiers and put $K_\eps=K*\rho_\eps$. Then $K_\eps\to K$ 
uniformly on compact sets of $\R^2$. 

Fix $u, \f\in\Hper$. Since $\Hper$ is continuously embedded into $C^0(\S^1,\R^2)$, then using $i)$ in Lemma \ref{L:2.5}
and (\ref{eq:dA}) with $K$ replaced by $K_\eps$, 
we have
\[
\begin{aligned}
A_K(u + \f) - A_K(u) &= A_{K_\e}(u + \f) - A_{K_\e}(u) + o(\e)= \int\limits_0^1 dA_{K_\e}(u+ s\f)\f\ ds+o(\eps) \\
&= \int\limits_0^1 ds\fint\limits_{\S^1} K_\e(u+s\f)\f\cdot(u' + s\f')~\! d\theta+o(\eps)\,.
\end{aligned}
\]
Taking the limit as $\eps\to 0$ we arrive at the identity
\[
A_K(u + \f) - A_K(u) = \int\limits_0^1 ds \fint\limits_{\S^1} K(u+s\f)\f\cdot(u' + s\f') ~\! d\theta \quad \text{ for any }u,\f \in \Hper\,.
\]
Since  $K$ is locally bounded, we deduce that 
\[
\big| A_K(u + \f) - A_K(u) - \irn K(u)\f\cdot iu'~\!d\theta \big|  = o (\|\f\|_{\Hper})\,.
\]
This implies that $A_K$ is Fr\'echet differentiable at $u$, with $dA_K(u)=\frac{1}{2\pi}K(u)iu'$. Since $K$ is continuous,
we also have that the function $u\mapsto K(u)iu'$ is continuous as a function $\Hper\to L^2(\S^1,\R^2)$, hence it
is continuous  $\Hper\to \Hdual$ as well. This concludes the proof of the first part of the lemma.

\medskip

Next, let  $K\in C^\infty_c(\R^2)$. By the remarks at the beginning of this section 
and by the linear dependence of the area functional from the weight function, we have
$$
dA_K(u)u=\irn K(u)u\cdot iu'~\!d\theta = -A_{\ddiv^*\!K}(u)= A_{\ddiv K}(u) + 2 A_K(u)
$$
for any $u\in \Hper$ (recall that  $\div (K(z)z)=-\ddiv^* K(z)=\ddiv K(z)+2K(z)$).

To conclude the proof for $K\in C^0(\R^2)\cap \widehat{\mathcal D}^1(\R^2)$ use 
the density result in Lemma \ref{L:dense} and $i)$ in Lemma \ref{L:2.5}.
\QED

The next Lemma evidently holds for curvatures $\rmH$ satisfying the assumptions in Theorem \ref{T:main}.

\begin{Lemma}
\label{L:iso_K}
Let $\rmH\in C^0(\R^2)$ be non constant. Assume that $(H_2)$ and $N_\rmH<\infty$ hold. Let $u\in\Hper\setminus\R^2$ and 
put $\rho_u=2\pi L(u)$ as in (\ref{eq:Linfty}). Then
\begin{gather}
\label{eq:AK_stima_vera}
\sqrt{4\pi}|A_{\rmH-1}(u)|\le \big(\int\limits_{\D_{\rho_{u}}(\p)}|\rmH(z)-1|^2~\!dz\big)^\frac12L(u)~,
\qquad 
|A_{\rmH-1}(u)|< N_\rmH~\!L(u)~\!,\\
\label{eq:estimate1_2}
|2A_\rmH(u)-dA_\rmH(u)u|\le N_\rmH~\!L(u)~\!.
\end{gather}
\end{Lemma}

\proof
By Lemma \ref{L:HinD} we have that $\rmH-1\in  C^0(\R^2)\cap \widehat{\mathcal D}^1(\R^2)\subset L^2(\R^2)$. 
Thus the inequalities in (\ref{eq:AK_stima_vera}) follow by using $ii)$ in Lemma \ref{L:2.5} and the Hardy inequality (\ref{eq:Hardy0}), compare with Lemma \ref{L:HinD_new}.

Further, Lemma \ref{L:2.6} gives that $A_\rmH$, $A_{\rmH-1}$ are differentiable, and 
$$
dA_{\rmH}(u)u=dA_1(u)u+dA_{\rmH-1}(u)u= 2A_1(u)+ A_{\ddiv\rmH}(u) + 2 A_{\rmH-1}(u)=
 A_{\ddiv\rmH}(u) + 2 A_{\rmH}(u)~\!.
$$ 
Thus $|2A_\rmH(u)-dA_\rmH(u)u|=|A_{\ddiv\rmH}(u)|$, so that (\ref{eq:estimate1_2}) follows 
from the weighted isoperimentric inequality (\ref{eq:isop}) with $\ddiv\rmH$ instead of $K$.
\QED

\section{The energy functional and proof of Theorem \ref{T:main}}
\label{S:energy}

Let $\rmH$ be a given continuous function on $\R^2$.  The energy functional
$$
E_{\rm H}(u)=L(u)+{A_\rmH}(u)\,,\qquad u\in { H}^1_{\rm per}
$$
is continuous on $\Hper$ and  {continuously} Fr\'echet differentiable on  ${\Hper}\setminus \R^2$. {Its differential  is} given by
$$
dE_\rmH(u)\f= \frac{1}{L(u)}\Big(\irn u'\f'~\!d\theta+L(u)\irn\rmH(u)\f\cdot iu'~\!d\theta\Big)\,,
$$
{use Lemma \ref{L:2.6} with $K$ replaced by $\rmH$}. If $u\in \Hper\setminus\R^2$ is a critical point for $E_\rmH$, then $u$ is a weak solution
to the system
$$
u''=L(u)\rmH(u) iu'.
$$
It easily follows that $|u'|$ is a constant. Precisely, $|u'|=L(u)$. Thus $u$ solves (\ref{eq:problem}), hence it is a 
$\rmH$-loop.

\medskip

Before going further, let us notice that $dE_{\rmH}(u)u=L(u)+dA_\rmH(u)u$ which, together with 
(\ref{eq:estimate1_2}), implies the crucial {estimate}
\begin{equation}
\label{eq:magic}
2E_{\rmH}(u)-dE_{\rmH}(u)u 
=L(u)+ 2A_\rmH(u)-dA_\rmH(u)u \ge (1-N_\rmH)L(u),
\end{equation}
which hold for any $u\in \Hper\setminus \R^2$.

Recall that a Palais-Smale sequence $u_n$ for $E_\rmH$ at a given energy level $c$, $(PS)_c$ sequence in brief, 
satisfies $u_n\in \Hper\setminus\R^2$,
$E_\rmH(u_n)=c+o(1)$ and $dE_\rmH(u_n)=o(1)$.

For completeness, we provide 
below the description of the behaviour of $(PS)_c$ sequences under the hypotheses $(H_1)$ and
$(H_2)$, including some details that are not  needed in the proof of Theorem \ref{T:main}. We start with the easiest case $\rmH\equiv 1$.

\begin{Lemma}
\label{L:H=1}
Let $c\in\R$ and let $u_n$ be a $(PS)_c$ sequence for $E_1$, such that $\overline u_n = 0$, see (\ref{eq:mean}).
Then there exist a subsequence $u_n$,  an integer $\ell\ge 1$ and $\theta_0\in\R$ such that $2c=\ell>0$ and
$u_n(\theta)\to e^{i(\theta-\theta_0)\ell}~~\text{in $H^1_{\rm per}$\,.}$
\end{Lemma}

\proof
The sequence $u_n/L(u_n)$ is bounded in $\Hper$. Thus 
$$
o(L(u_n))=dE_1(u_n)u_n=L(u_n)+\irn u_n\cdot iu'_n~\!d\theta \ge L(u_n)-\irn |u_n| |u'_n|~\!d\theta\ge L(u_n)-2\pi L(u_n)^2,
$$
by (\ref{eq:Linfty}) and by the Cauchy-Schwarz inequality.
We infer that $L(u_n)$ can not converge to zero. 

In addition, we notice that 
$$
2c + o(1) = 2E_1(u_n) = L(u_n) + dE_1(u_n)u_n =L(u_n) (1 + o(1))\,.
$$
Thus $L(u_n)\to 2c>0$,  the sequence $u_n$ is bounded in $H^1_{\rm per}$ and we can assume that
$u_n\to U$ weakly in $H^1_{\rm per}$. In fact, $u_n\to U$  in the $H^1_{\rm per}$- norm, because
$$
o(1)=dE_1(u_n)(u_n-U)= \frac{1}{L(u_n)}\irn u'_n(u_n-U)'~\!d\theta+\irn (u_n-U)\cdot iu'_n~\!d\theta=
\frac{L(u_n- U)^2}{L(u_n)}+o(1)\,.
$$
The strong convergence gives $L(U)=2c$, thus $U$ is not constant. Since $E_1$ is of class $C^1$ in 
$\Hper\setminus\R^2$, we see that $dE_1(U)=0$. Hence $U\in\Hper$ is a non constant
solution to $U''=L(U) iU'$. The conclusion of the proof follows via Fourier expansion. 
\QED

Next, we deal with non constant curvatures.

\begin{Lemma}
\label{L:HPS}
Assume that $\rmH \in C^0(\R^2)$ satisfies $(H_1)$ and $(H_2)$. Let $c\in\R$ and let 
$u_n$ be a $(PS)_c$ sequence for $E_\rmH$. 
Then $c>0$ and there exist a subsequence $u_n$ which satisfies one of the next alternatives:
\begin{enumerate}
\item[a)] 
the sequence of means $\p_n\subset \R^2$ is unbounded, and 
$u_n-\p_n$ converges in $\Hper$ to a  parametrization of the unit circle about the origin of topological degree $\ell\ge 1$.
In particular, $c=\frac{\ell}{2}\ge \frac12$;
\item[b)] 
$u_n\to U$ in $H^1_{\rm per}$, where $U$ is a  $\rmH$-loop. 
\end{enumerate}
\end{Lemma}

\proof
We start by noticing the crucial inequality 
\begin{gather}
\label{eq:magic_n}
2c-dE_{\rmH}(u_n)u_n
\ge (1-N_\rmH)L(u_n) +o(1),
\end{gather}
compare with  (\ref{eq:magic}). 
By adapting an 
argument already used in the proof of Lemma \ref{L:H=1}, we  show  that 
\begin{equation}
\label{eq:liminf}
\liminf_{n\to\infty} L(u_n)>0.
\end{equation}
In fact, the sequence $\tfrac{u_n-\p_n}{L(u_n)}$ is bounded in ${\Hper}$. Thus 
$$
o(L(u_n)) = dE_\rmH(u_n)(u_n - \p_n) \ge L(u_n) - \irn |\rmH(u_n)||u_n - \p_n||u_n'| d\theta.
$$
The function $\rmH$ is bounded, since it is continuous and satisfies $(H_2)$. Since 
$\|u_n-\p_n\|_\infty\le 2\pi L(u_n)$,  we infer the estimate
$o(L(u_n))\ge L(u_n)-2\pi\|\rmH\|_\infty L(u_n)^2$,
which ends the proof of (\ref{eq:liminf}).

\medskip 

We divide the rest of the proof in three steps.

\paragraph{Step 1:}{\em The sequence $L(u_n)$ is bounded.}~\\
We have that $dE_{\rmH}(u_n)(u_n-\p_n)=o(L(u_n))$ because the sequence $\tfrac{u_n-\p_n}{L(u_n)}$ is bounded in ${\Hper}$. Assume by contradiction that, for a subsequence,  $L(u_n)\to \infty$. Then
(\ref{eq:magic_n}) easily implies
\begin{equation}
\label{eq:contradicts}
 - dE_{\rmH}(u_n)\p_n \ge (1-N_\rmH)L(u_n)+o(L(u_n)).
\end{equation}
It follows that the sequence $\p_n/L(u_n)$ can not be bounded in $\R^2\subset\Hper$, as
$dE_\rmH(u_n)=o(1)$. Hence 
\begin{equation}
\label{eq:pL}
L(u_n)=o(|\p_n|),
\end{equation}
and in particular $|\p_n|\to \infty$. By the triangle inequality and (\ref{eq:Linfty}), we have 
$$
|u_n|\ge |\p_n|-|u_n-\p_n|\ge |\p_n|\Big(1-\frac{2\pi L(u_n)}{|\p_n|}\Big)\qquad\text{on $\S^1$,}
$$
which, together with (\ref{eq:pL}), implies that $2|u_n|\ge |\p_n|$ for $n$ large enough. We also infer that
$|u_n|\to\infty$ uniformly on $\S^1$.

Trivially, $dE_1(u_n)=dE_{\rmH}(u_n)-dA_{\rmH-1}(u_n)$ vanishes on constant functions. Thus, 
we can  estimate
\begin{equation}
\begin{aligned}
\label{eq:gives}
\big|dE_{\rmH}(u_n)\p_n\big|&= \big|\irn (\rmH(u_n)-1)\p_n\cdot iu'_n~\!d\theta\big|\le 
2\irn |u_n|~\!|\rmH(u_n)-1||u_n'|~\!d\theta \\
&\le 2L(u_n)\Big(\irn |u_n|^{2}|\rmH(u_n)-1|^2 ~\!d\theta\Big)^\frac12~\!.
\end{aligned}
\end{equation}
The last integral in (\ref{eq:gives}) converges to zero by assumption $(H_2)$. Therefore 
$dE_{\rmH}(u_n)\p_n=o(L(u_n))$, which contradicts (\ref{eq:contradicts}) because $N_\rmH<1$, and concludes Step 1.

\paragraph{Step 2:} {\em If $\p_n$ is unbounded, then, 
up to a subsequence, the alternative $a)$ occurs.}~\\
By Step 1, we can assume that the sequence $u_n-\p_n$ converges weakly in $\Hper$. Since
$u_n-\p_n$ converges uniformly on $\S^1$, using also (\ref{eq:AK_stima_vera}) we see that there exists $R>0$ such that 
$$
\sqrt{4\pi}|A_{\rmH-1}(u_n)|\le L(u_n)\big(\int\limits_{\D_{R}( \p_n)}|\rmH(z)-1|^2~\!dz\big)^\frac12.
$$
We infer that  $A_{\rmH-1}(u_n)=o(1)$, because $|\p_n|\to\infty$ and $\rmH-1\in L^2(\R^2)$ by Lemma \ref{L:HinD_new}.
Therefore
$$
E_1(u_n-\p_n)=E_1(u_n)=E_\rmH(u_n)-A_{\rmH-1}(u_n)= c + o(1)\,.
$$
In addition, for any $v\in H^1_{\rm per}$ we can estimate
$$
|dA_{\rmH-1}(u_n)v| =\Big|\irn (\rmH(u_n)-1)v\cdot iu'_n~\!d\theta\Big|\le \|\rmH\circ u_n-1\|_\infty L(u_n)\|v\|_2 =o(\|v\|_{\Hper}),
$$
by $(H_2)$ and since $|u_n|\to\infty$ uniformly. Thus $dA_{\rmH-1}(u_n)=o(1)$ in $\Hdual$, which implies
$$
dE_1(u_n-\p_n)=dE_1(u_n)=dE_\rmH(u_n)-dA_{\rmH-1}(u_n)=o(1)\quad\text{ in $\Hdual$.}
$$
We showed that $u_n-\p_n$ is a $(PS)_c$ sequence for $E_1$, which concludes Step 2, thanks to Lemma \ref{L:H=1}.

\paragraph{Step 3:}
{\em If $\p_n$ is bounded, then, 
up to a subsequence,
the alternative $b)$ occurs.}~\\
By Step 1 we can assume that $u_n\to U$ weakly in $\Hper$. 
Thus
$$
\begin{aligned}
o(1)=dE_\rmH(u_n)(u_n-U)&=\frac{1}{L(u_n)}\irn u'_n\cdot (u_n-U)'~\!d\theta+dA_\rmH(u_n)(u_n-U)=\frac{L(u_n-U)^2}{L(u_n)}+o(1)
\end{aligned}
$$
because $2\pi dA_\rmH(u_n)=\rmH(u_n)iu'_n$ is bounded in $L^2(\S^1,\R^2)$ and $u_n-U\to 0$ in $L^2(\S^1,\R^2)$.
We infer that $u_n\to U$ strongly in $H^1_{\rm per}$, and thus $L(u_n)=L(U)+o(1)$. It follows that 
$U$ is non constant by (\ref{eq:liminf}). Then, by continuity we also have 
$E_\rmH(U)=c$ and
$dE_\rmH(U)=0$, that is, $U$ is an $\rm H$-loop. This ends Step 3.

\medskip

Finally, we notice that $c\ge \frac12>0$ if the first alternative occurs. Otherwise, let $U$ be the loop in $b)$.  Then 
(\ref{eq:magic}) gives
$2c=2E_\rmH(U)-dE_{\rmH}(U)U \ge (1-N_\rmH)L(U)>0$, which implies $c>0$.
 The lemma is completely proved. \QED

Before proving Theorem \ref{T:main} we point out a lemma about regular parametrizations of circles in $\R^2$. 

\begin{Lemma}
\label{L:negative}
Assume that $\rmH \in C^0(\R^2)$ satisfies $(H_1)$ and $(H_2)$ and let $R\geq 2(1+N_\rmH)$. If $\rmH$ is non constant then 
$E_\rmH(Re^{i\theta}+p)<0$ for any $p\in\R^2$.
\end{Lemma}

\proof
Let us start with some computations which hold for any $R>0$. 
The loop $\omega(\theta)=Re^{i\theta}+p$ parametrizes $\partial\D_R(p)$, has constant scalar speed $|\omega'|=R$ and
evidently $-j_{\omega}$ is the characteristic function of $\D_{R}(p_\eps)$. Thus we can compute
\begin{equation}
\label{eq:utile}
L(Re^{i\theta}+p)=R~,\qquad 2\pi A_K(Re^{i\theta}+p)=-\int\limits_{\D_{R}(p_\eps)}K(z)~\!dz
\quad\text{for any $K\in L^2_{\rm loc}(\R^2)$. }
\end{equation}
Using the first equality in (\ref{eq:utile}) and (\ref{eq:AK_stima_vera}) we get
$
|A_{\rmH-1}(Re^{i\theta}+p)|< RN_\rmH$. This allows us to estimate
$$
E_\rmH(R{e^{i\theta}}+p)=E_1(R{e^{i\theta}})+
A_{\rmH-1}(R{e^{i\theta}}+p)<(1+N_\rmH)R-\frac{R^2}{2}.
$$
The conclusion for $R\geq 2(1+N_\rmH)$ readily follows. 
\QED

\subsection{Proof of Theorem \ref{T:main}}
\label{SS:Proof}

We can assume that $\rmH$ is non constant, otherwise any circle of radius $1$ is an $\rmH$-loop.

Let $\tilde R =  2(1+N_\rmH)$. We define 
\begin{equation}
\label{eq:mp}
\Gamma=\big\{\gamma\in C^0([0,1], {\Hper})~|~ \gamma(0)\in \R^2~,\quad \gamma(1)(\theta)=\tilde R{e^{i\theta}} \}~,
\quad \cmp=\inf_{\gamma\in\Gamma}\max_{t\in[0,1]} E_\rmH(\gamma(t)).
\end{equation}
We divide the proof in few steps.

\paragraph{Step 1:} {\em If $(H_1)$ and $(H_2)$ hold, then there exists a Palais-Smale sequence at the level 
$\cmp$.}~\\
Constant functions are local minima for $E_\rmH$.
This is a consequence of the isoperimetric inequalities in (\ref{eq:AK_stima_vera}), (\ref{eq:iso}),  which give
$$
E_{\rmH}(u) =  \big[L(u)+A_{\rmH-1}(u)\big]+A_1(u)
\ge(1-N_\rmH)  L(u)- \frac12 L(u)^2 \quad \text{for any $u\in\Hper$}.
$$
Further, $L(\tilde Re^{i\theta})= \tilde R >1-N_\rmH$ by (\ref{eq:utile}).
Thus, any path $\gamma\in \Gamma$ crosses $\{L(u)=1-N_\rmH\}$, which implies  $\cmp>0$. Since in addition
$E_{\rmH}(\tilde Re^{i\theta})<0$ by Lemma \ref{L:negative}, we see that $\cmp$ is a mountain pass level for $E_\rmH$.

A deformation lemma for $C^1$ functionals based on pseudo gradient vector fields, see for instance \cite[Theorem 8.2]{AM},
provides the existence of a $(PS)_{\cmp}$ sequence for $E_\rmH$.

\paragraph{Step 2:} {\em If $(H_1)$ and $(H_2)$ hold, then $\cmp \leq \frac{1}{2}$.}~\\
Fix $\e>0$ and use $(H_2)$ to find $p_\e \in \R^2$ such that
\[
|p_\e| > 2\tilde R + \frac{\tilde R^2}{\e}~\!C_\rmH~,\qquad C_\rmH:=\sup_{z\in\R^2}|(\rmH(z)-1)z|
\,.
\]
Consider the path $\gamma_\e \in \Gamma$ given by
\[
\gamma_\e(t)(\theta) = 
\begin{cases}
2t\tilde R{e^{i\theta}}+p_\e  & t\in \big[0, \frac{1}{2}\big)\\
\tilde R{e^{i\theta}}+2(1-t)p_\e  & t\in \big[ \frac{1}{2}, 1\big]\,.
\end{cases}
\]
If $t \in (\frac{1}{2}, 1]$ then $E_{\rmH}(\gamma_\e(t) ) < 0$  by Lemma \ref{L:negative}.
If $t \in (0, \frac{1}{2}]$, then $\gamma_\eps(t)$ parametrizes the circle of radius $2t\tilde R$ about $p_\eps$.
Since $|p_\eps|>2\tilde R$, we have that $\D_{2t\tilde R}(p_\e)\subset \D_{\tilde R}(p_\e)\subset \{|z|>|p_\eps|/2\}$. Thus
$|\rmH(z)-1|\le C_\rmH|z|^{-1}\le 2C_\rmH|p_\eps|^{-1}$ on $\D_{2t\tilde R}(p_\e)$ and therefore
$$
2\pi |A_{\rmH-1}(\gamma_\eps(t))| =
\big| \int\limits_{\D_{2t\tilde R}(p_\e)}(\rmH(z)-1)~\!dz\big| \le 
  \int\limits_{\D_{\tilde R}(p_\e)}|\rmH(z)-1|~\!dz\le 2\pi \frac{ C_\rmH \tilde R^2}{|p_\eps|}<2\pi~\!\eps
$$ 
by (\ref{eq:utile}). It follows that 
$
E_{\rmH}(\gamma_\e(t) ) = L(\gamma_\eps(t))+A_1(\gamma_\eps(t))+A_{\rmH-1}(\gamma_\eps(t))
\le 2t\tilde R-2t^2\tilde R^2 +\eps  \leq \frac{1}{2}+\e\,,
$
and we can conclude that 
$$
\cmp\le \sup_{t\in[0,\frac12]} E_{\rmH}(\gamma_\eps(t) ) \le  \frac{1}{2}+\e.
$$
Since $\eps$ was arbitrarily chosen, this proves that  $\cmp \leq \frac12$, as claimed. 

\paragraph{Step 3:}{\em If $(H_1)$, $(H_2)$ and $(H_3)$ hold, there exists an $\rmH$-loop with energy $\cmp$.}~\\
Step 1 provides the existence of a $(PS)_{\cmp}$ sequence.
If $\cmp<\frac12$, then the existence of a non constant $\rmH$-loop is given by Lemma \ref{L:HPS}.

Otherwise, $\cmp=\frac12$ by Step 2. 
To conclude the proof we  show that,  in this case, the last assumption $(H_3)$ implies that 
there exists a circle of radius $1$ which is also an $\rm H$-loop.

Let $\tilde{p}\in \R^2$ be given by $(H_3)$, so that $\rmH(z)\ge  1$ on $\D_{\tilde R}(\tilde{p})$. Consider the path $\gamma \in \Gamma$ given by
\[
\gamma(t) = 
\begin{cases}
2t\tilde{ R}{e^{i\theta}}+\tilde{p}  & t\in \big[0, \frac{1}{2}\big)\\
\tilde{ R}{e^{i\theta}}+2(1-t)\tilde{p} & t\in \big[ \frac{1}{2}, 1\big]\,.
\end{cases}
\]
Since $E_{\rmH}(\gamma(t) )<0 = E_{\rmH}(\gamma(0) )$ for any $t \in [\frac12, 1]$ by Lemma \ref{L:negative}, we see that
there exists $t_0 \in (0, \frac{1}{2})$ such that $\max_{t \in [0,1]} E_{\rmH}(\gamma(t)) = E_{\rmH}(\gamma(t_0) )$.
Using the identity $E_\rmH= E_1+A_{\rmH-1}$ and (\ref{eq:utile})  (recall that $\rmH-1\in L^2(\R^2)$ by Lemma 
\ref{L:HinD_new}), we can compute 
\[
\frac12 = \cmp \leq E_{\rmH}(\gamma(t_0) ) = 2t_0{\tilde R}-2t^2_0{\tilde R}^2 - \frac{1}{2\pi}\int\limits_{\D_{2t_0{\tilde R}}(\tilde p)}(\rmH(z)-1)~\!dz \leq (2t_0{\tilde R})-\frac12(2t_0{\tilde R})^2  \le \frac{1}{2}\,.
\]
Thus equalities hold everywhere in the  above formula. In particular, we infer that 
$2t_0 {\tilde R}= 1$ and that the continuous function $\rmH - 1$ vanishes on $\partial\D_1(\tilde{p})$. 
Therefore, the loop $u(\theta)=\tilde{p}+ e^{i\theta}$ is an $\rm H$-loop. 

\medskip

Theorem \ref{T:main} is completely proved. \QED

\appendix

\section{\hskip-0.6cm ppendix. Final remarks}

Theorem \ref{T:main} has the next straightforward extension. 

\begin{Corollary}
\label{C:Hinfty}
Assume that $\rmH\in C^0(\R^2)$ satisfies $(H_1)$ and  

$(H^\lambda_2)$ $\rmH(z)-\lambda=o(|z|^{-1})$ as $|z|\to\infty$, for some $\lambda\in \R, \lambda\neq 0$;

$(H^\lambda_3)$ $\lambda^{-1}\rmH(z)\ge  1$ if $|\lambda||z-\tilde{p}|< 2(1 + N_\rmH)$.

Then there exists at least one $\rmH$-loop.
\end{Corollary}

\proof
Recall that changing the orientation of a  curve  
changes the sign of its curvature. Thus we can assume $\lambda >0$. Since the function
$$
\rmH_\lambda(z)=\frac{1}{\lambda}\rmH\big(\frac{z}{\lambda}\big)\,,
$$
satisfies $N_{\rmH_\lambda}=N_\rmH<1$, $(H_2)$ and $(H_3)$ (with $\lambda\tilde p$ instead of $\tilde p$), then
Theorem \ref{T:main} gives the existence of a $\rmH_\lambda$-loop $u$.
To conclude the proof, it
suffices to check that $u_\lambda:=\lambda^{-1} u\in\Hper$ is a solution to 
$u_\lambda''=L(u_\lambda)\rmH(u_\lambda)iu_\lambda'$.
\QED

The assumption $\lambda\neq 0$ in Corollary \ref{C:Hinfty} is needed, because of the next nonexistence result.

\begin{Theorem}
\label{T:ne}
Let $\rmH\in C^0(\R^2)\cap L^2(\R^2)$. If $(H_1)$ holds, then  no $\rmH$-loop exist.
\end{Theorem}

\proof
We have that $\rmH\in \widehat{\mathcal D}^1(\R^2)$ and  $\|\ddiv^*\rmH\|_2=\|\ddiv\rmH\|_2=N_\rmH$ by $ii)$ in Lemma \ref{L:dense}. 

If $u\in \Hper$ solves
$u''=L(u){{\rm \rmH}}(u)iu'$, then
$$
L(u)^2=\irn|u'|^2 d\theta =-L(u)\irn \rmH(u)u\cdot iu' d\theta = L(u) A_{\ddiv^*\!\rmH}(u)\le  \frac{1}{\sqrt{4\pi}} \|\ddiv^*\rmH\|_2 L(u)^2 = N_\rmH L(u)^2
$$
by (\ref{eq:isop}).
This implies that $u$ is constant, because $N_\rmH<1$.
\QED

We now provide few remarks and examples to comment our main hypotheses (the quantity $M_\rmH$ is defined
in (\ref{eq:MH})).

\begin{Claim}
\label{E:Ex1}
The assumptions $M_{\rmH}<1$ and $N_{\rmH}<1$ are not comparable.
\end{Claim}

\proof Let $\beta>1$, $t>0$. The curvature
$$
{\rmH_{\beta,t}}(z)=\begin{cases}1+t|z|^{-\beta}&\text{if $|z|\ge {1}$}\\
1+t(2-|z|^\beta)&\text{if $|z|< {1}$}
\end{cases}
$$
{is of class $C^1$, satisfies $(H_2)$, $(H_3)$ and }
$$
M_{\rmH_{\beta,t}}=\beta t~,\qquad N_{\rmH_{\beta,t}}=M_{\rmH_{\beta,t}}~\!\sqrt{{\frac{\beta}{2(\beta^2-1)}} }. 
$$
In particular, 
$$
\begin{aligned}
N_{\rmH_{\beta,t}}<1\leq M_{\rmH_{\beta,t}} \quad&\text{if $\beta$ is large enough and $\beta{\leq}t^2\beta^3<2(\beta^2-1)$}\\
M_{\rmH_{\beta,t}}<1\leq N_{\rmH_{\beta,t}} \quad&\text{if $\beta$ is close to $1^+$ and $2(\beta^2-1)\leq t^2\beta^3< \beta$.}
\end{aligned}
$$

\vspace{-0.7cm}
\QED

\bigskip

The curvature ${\rmH_{\beta,t}}$ in the previous example is radially symmetric and evidently there exists $r_{\beta,t}\in(0,1)$
such that the circle of radius $r_{\beta,t}$ about the origin can be parametrized by a ${\rmH_{\beta,t}}$-loop. To exhibit 
examples of curvatures satisfying $N_\rmH<1$ and/or  $M_\rmH<1$ for which the existence of $\rmH$-loops is
not evident, one can consider curvatures of the type
$\rmH_{\beta,t}+\eps\f$, where $\f$ is a generic function in $C^\infty_c(\R^2)$ and $\eps>0$.

\medskip

\begin{Claim}
\label{E:Ex2}
For any $\delta>0$, there exists a continuous curvature $\rmH$ satisfying $1<N_{\rmH}<1+\delta$ and $(H_2)$, 
for which $E_\rmH$ admits a Palais-Smale sequence $u_n\in\Hper $ such that $L(u_n)\to \infty$.
\end{Claim}

\proof
For any $\eps\in (0,1)$ we introduce the radial curvature
$$
{\rmH}_\eps(z)=1+\frac{6}{3-\eps^2}~\!\psi_\eps(r)~,
\quad\text{where}\quad
\psi_\eps(r)=
\begin{cases} 1+|\log \eps|-\frac{r}{\eps}&\text{if $r\le\eps$}\\
|\log r|&\text{if $\eps<r\le1$}\\
0&\text{if $r>1$}
\end{cases}.
$$
Assumption $(H_2)$ is trivially satisfied. A simple computation gives 
$$
{N_{\rmH_\eps}^2}= \frac92~\!\frac{2-\eps^2}{(3-\eps^2)^2}\searrow 1\qquad \text{as $\eps\searrow 0$.}
$$
The curve $u_n(\theta)=e^{in\theta}$ parameterizes the unit circle and has degree $n$. In fact
$u_n$ is a 
$\rmH_\eps$-loop, that is, $dE_{\rmH_\eps}(u_n)=0$, because $\rmH_\eps(u_n)\equiv 1$. Since
in addition 
$$
E_{\rmH_\eps}(u_n)=n{E_{\rmH_\eps}(u_1)}= n\big(1-\frac{1}{2\pi}\int\limits_{\D_1} \rmH_\eps(z)~\!dz\big)
= n\big(\frac12-\frac{6}{3-\eps^2}\int\limits_{0}^1 r\psi_\eps(r)~\!dr\big)=0,
$$
it turns out that $u_n$ is a Palais-Smale sequence for $E_\rmH$. However, $L(u_n)=n\to\infty$.
\QED

\end{document}